\documentclass[runningheads,a4paper]{llncs}

\usepackage{amssymb}
\usepackage{enumitem}
\usepackage{graphicx}
\usepackage{amsmath}
\usepackage[numbers]{natbib}
\usepackage{mathtools}
\usepackage{subfigure}
\usepackage{comment}
\usepackage[linesnumbered,ruled,vlined]{algorithm2e}
\usepackage{url}
\usepackage{bbm}

\newcommand{\keywords}[1]{\par\addvspace\baselineskip
\noindent\keywordname\enspace\ignorespaces#1}

\begin{document}

\mainmatter  
\title{Understanding Spectral Graph Neural Network}
\author{Xinye Chen%
\thanks{xinye.chen@manchester.ac.uk}}

\institute{Department of Mathematics, University of Manchester\\
Manchester, M13 9PL, United Kingdom
}

\maketitle

\begin{abstract}
Graph neural networks have developed by leaps and bounds in recent years due to the restriction of traditional convolutional filters on non-Euclidean structured data. Spectral graph theory mainly studies fundamental graph properties using algebraic methods to analyze the spectrum of the adjacency matrix or Laplacian matrix of a graph, which lays the foundation of graph convolutional neural networks. This report is more than notes and self-contained which comes from my Ph.D. first-year report literature review part, it illustrates how the graph convolutional neural network model is motivated by spectral graph theory, and discusses the major spectral-based models associated with their fundamentals. The practical applications of the graph convolutional neural networks defined in the spectral domain are also reviewed.
\end{abstract}

\section{Overview}\label{sec:intro}

In recent years, the growing computing power of machines has been greatly accelerating the development of deep learning. With the advance of computational hardware and research output, deep learning has made great progress and achieved great success in many fields including translation, object recognition, recommendation systems, and so on. Convolutional Neural Networks (CNNs) are effective deep learning techniques for addressing numerous machine learning and data mining problems, achieving promising performance in image processing \cite{6795724, NIPS2012_4824}, document recognition \citep{726791, kim-2014-convolutional}, object recognition \citep{9050545, 7780459, NIPS2015_14bfa6bb}, speech recognition \citep{6296526, DBLP:conf/ssw/OordDZSVGKSK16}, game of Go \citep{articlesdha}, and bioinformatics \citep{10.3389/fgene.2021.690049}, in which the data are associated with an underlying grid-like structure. We categorize such data with grid-like structure into the class of Euclidean data, in which we can operate the computation with standard inner products, subtract one vector from another, apply matrices to vectors, etc. For example, data like time signals and images{---}that are discretized on regular Cartesian grids{---}can be applied to operations like convolution by simply sliding the same window over the signal and computing inner products. 

Euclidean data is often easy to be manipulated by networks with convolutional architectures \citep{6795724cnn} because of the \emph{translational equivariance} and \emph{invariance} properties arising from such grid structure \cite{af52b944fc104958bdc3e16d4ee15903}. However, the nature of data defined on non-Euclidean domains like graphs and manifolds indicates that there are no such familiar properties as global parameterization, a common system of coordinates, vector space structure, or shift-invariance  (we refer the reader to \cite{DBLP:journals/corr/abs-1904-11486,DBLP:journals/corr/WorrallGTB16} for further details), e.g., the characteristics of users in social networks can be modeled as signals on the vertices of the social graph \citep{LDPA}; papers linked to each other via citations can be categorized into different groups according to topics  \citep{Velickovic2018GraphAN}; traffic data of different roads and times can be modeled as graph structure signals \citep{TGCNNADF}. Therefore, the practice of deep learning on Euclidean data remains a popular topic waiting for optimal solutions.

A neural network structure that can efficiently operate and extract useful features on such non-Euclidean-domain data is very desired. 

\subsection{Motivation}\label{subsec:motivation}
Broadly speaking, graphs are ubiquitous in the real world in the form of representing objects associated with their relationships such as social networks, e-commerce networks, biology networks, and traffic networks (see \cite{Zhang2018DeepLO} for a review). Technically, graphs are generic data representation forms that are useful for illustrating the geometric structures of data domains in a great number of applications, including social, energy, transportation, sensor, and neural networks \citep{cite2}. The graph neural network is motivated by CNNs which have been successfully applied in the field of computer vision \citep{5537907cnn, 726791, GNNAMA, 6338939objectdetection}. CNNs are essentially a high-performance end-to-end learning framework \footnote{End-to-end learning refers to training a possibly complex learning system by applying gradient-based learning to the system as a whole \citep{Glasmachers2017LimitsOE}.} for processing image information, but it can only operate on regular Euclidean data like 2D grid and 1D sequence \citep{GNNAMA}. Besides, the characteristics of CNNs: local connection, shared weights and the use of multi-layer are of great importance in addressing problems in graph domain \citep{DLYANN}, \citep{GNNAMA}, because 1) graphs are the most typical locally connected structure; 2) shared weights reduce the computational cost compared with traditional spectral graph theory \citep{SPinbook, GNNAMA}. To introduce the graph neural network, we need first to associate it with spectral graph theory, whose focus is to examine the eigenvalues (or spectrum) of a matrix (usually Laplacian matrix) associated with a graph and utilize them to determine the structural properties of the graph \citep{SPinbook}.

Graph deep learning (or geometric deep learning) is a hyperonym for emerging techniques attempting to generalize (structured) deep neural models to non-Euclidean domains such as graphs and manifolds \citep{af52b944fc104958bdc3e16d4ee15903}. The notation of graph neural networks was first mentioned in \cite{1555942}, and further developed and completed in \cite{4700287}. These early works presented graph neural networks that need computationally expensive training that learn the target node's representation by propagating neighbor vertex or link information via recurrent neural networks in an iterative way until a stable convergence is achieved \citep{dai2018learning, li2016gated, Wu2021ACS}. 

Currently, most deep learning methods such as LSTM and CNN are good at processing sequence data, image data, video data, text data, and others defined in the Euclidean domain. However, most deep learning algorithms do not perform very well with data on non-Euclidean domains. By contrast, graph neural networks, which are the current popular topic in the deep learning area, can achieve a good performance on non-Euclidean domains. 
 
\subsection{Related work on graph}\label{subsec:background}
Graph neural networks can be provided a taxonomy that divides graph neural networks into five categories, graph convolutional networks, graph attention networks, graph autoencoders, and graph generative networks \cite{Wu2021ACS}. Here, we mainly introduce graph convolutional neural networks on the spectral domain as well as the basics of spectral graph theory. 

The recent years of graph convolutional neural networks (GCNs) can be listed as the following; The first work on spectral GCNs can be traced back to \cite{ae482107de73461787258f805cf8f4ed} which is based upon a hierarchical clustering of the domain, and another based on the spectrum of the graph Laplacian respectively. Then, to avoid high computational complexity arising from eigen decomposition, Chebyshev GCN (ChebNet) utilizes truncated expansion of Chebyshev polynomials \citep{Hammond:131283} to fit convolution kernels \citep{10.5555/3157382.3157527}. In the meantime, another work \citep{Susnjara2015} proposes an accelerated algorithm based on the Lanczos method that adapts to the Laplacian spectrum without explicitly computing it and achieves higher accuracy without increasing the overall complexity significantly compared to methods based on Chebyshev polynomials. GCN, a scalable approach for semi-supervised learning on graph-structured data that is based on an efficient variant of convolutional neural networks, can operate directly on graphs \cite{Kipf:2016tc}. Graph Convolutional Recurrent Network (GCRN), a generalization of classical recurrent neural networks (RNN) and graph CNN, can predict structured sequences of data, which represent series of frames in videos, spatio-temporal measurements on a network of sensors, or random walks on a vocabulary graph for natural language modeling \citep{journals/corr/SeoDVB16}. CayleyNets GCN introduces a new spectral-domain convolutional architecture for deep learning on graphs based on Cayley filters instead of Chebyshev filters \citep{CayleyNets2017}. 

The disadvantage of spectral-based GCN is that the learned filters rely on the Laplacian eigenbasis, depending on the graph structure, which in turn means a model trained on a specific structure that can not be directly extended to another graph with a different structure \citep{Velickovic2018GraphAN}. Besides, early work on spectral GCNs is limited to undirected graphs. 

\subsection{Tasks on graph}\label{subsec:overview}

In practice, the graph structure itself can be categorized into homogenous or heterogeneous levels \citep{Zhang2018DeepLO, GNNBook2022}. A graph is heterogeneous if each node and each edge are associated with a type and there is more than one type that exists in the graph nodes or edges, otherwise, the graph is homogenous. A more formal definition can be referred to \cite{hamiltomgrl}. More categories with respect to graph structure can be referred to \cite{hamiltomgrl, ma2021deep, GNNBook2022}. 

On top of the difference in graph structure, GCN approaches can be classified into two categories, spectral-domain and spatial-domain methods \cite{hamiltomgrl}. Here, we introduce the basics of spectral graph theory according to \cite{Chung:1997} which is spectral GCNs based, and review the methods of graph GCN, with a focus on the spectral domain. With respect to all the graph application tasks, introduced here, we assume the graph is homogenous and the task is node-level. The paper is organized as follows. The first two sections review the motivation and background of GCNs while briefly discussing the categories in graph-related tasks.   

Besides, we can divide the graph-related tasks into three categories, namely node-level, edge-level, and graph-level, which allows us to focus on different graph analytics tasks \citep{Wu2021ACS}: 
\begin{enumerate}
\item \textit{Node-level}: tasks about predicting the node for regression or classification. In this task, the entire data are stored in a graph, and each node is an individual sample, e.g., node-level anomaly detection \citep{9565320};  
\item \textit{Edge-level}: tasks about predicting the edge or link for classification, e.g., relation prediction in knowledge graphs \citep{nathani-etal-2019-learning}, criminal intelligence analysis \citep{10.1371/journal.pone.0154244}, protein--protein interaction \citep{ijcai2021p506, Jha2022};
\item \textit{Graph-level}: tasks about predicting graph for classification. This task requires pooling techniques, e.g., introducing graph pooling layers, to obtain the representation of a whole graph. Graph level tasks include graph-level anomaly detection \citep{9565320, zhang2022dualdiscriminative, 10.1145/3580305.3599524}, neural machine translation \citep{beck-etal-2018-graph}, molecular property prediction \citep{Jiang2021}, etc.  
\end{enumerate}

\section{Basic graph concepts}\label{sec:basics}
A graph can be defined as $\mathcal{G} = (\mathcal{V},\mathcal{E},A)$, where $\mathcal{V}$ is the set of vertices or nodes, $\mathcal{E}$ is set of edges or links, and $A$ is the adjacency matrix of size $n \times n$. $v_{i} \in \mathcal{V}$ denotes a node, and $e_{i,j} \in \mathcal{E}$ denotes an edge connecting $v_{i}$ and $v_{j}$ in a graph $\mathcal{G}$. If an edge $e_{i,j}$ exists in graph, denoted by $e_{i,j} \in \mathcal{E}$, then $A_{i,j} > 0$,  otherwise $A_{i,j} = 0$ and $e_{i,j} \notin \mathcal{E}$.

\textbf{Degree of vertex:} The degree of node $i$ is $d_{i}$, representing the number of edges connected to node $i$,which is defined by
\begin{equation}
d_{i} = \sum_{j=1}^{n}\mathbbm{1}_{\mathcal{E}} \{e_{i,j}\},
\end{equation}
where $\mathbbm{1}$ is indicator function.

Given a graph $\mathcal{G}$, the degrees matrix $D \in \mathbb{R}^{n \times n}$ is
\begin{equation}
D_{i,j} = \left\{
\begin{array}{rcl}
d_{i} & & \textrm{if $i = j$}\\
0 & & \textrm{otherwise}\\
\end{array} \right..
\end{equation}

Undirected graph is graph with undirected edges and has $A_{i,j} = A_{j,i}$. In contrast, directed Graph is graph with directed edges, which may not satisfy $A_{i,j} \ne A_{j,i}$. The spectral graph convolutional network is defined on an undirected graph. In fact, an undirected graph is a special case of a directed graph.

\textbf{Diameter of graph \cite{weighteddiametergraph}: }  Given a connected graph $\mathcal{G}$, for two vertices $v_{a}$ and $v_{b} \in \mathcal{V}$, a path between $v_{a}$ and $v_{b}$ is a sequence $\pi = (e_{1}, e_{2},\ldots,e_{k})$ where $e_{i}=(v_{i-1}, v_{i}) \in \mathcal{E}$ and $v_{i} \in \mathcal{V}$ for $i \in {1,\ldots,k}$ with $v_{0}=v_{a}$ and $v_{k}=v_{b}$. We denote $e \in \pi$ if the edge $e \in \mathcal{E}$ belongs to the path $\pi$, i.e., if $e=e_{i}$ for an $i \in {1,\ldots,k}$. The distance between two vertices $v_{i}$ and $v_{j}$ is the number of edges in $\mathcal{E}$ in the shortest path connecting these two vertices, denoted by
\begin{equation}
dist(v_{i}, v_{j}) = \min\sum_{e \in \pi}w_{e},
\end{equation}
where $w_{e}$ is the weight on the edge $e$, $w_{e} = 1$ if it applies to unweighted graph.

The diameter of $\mathcal{G}$, denoted by $diam(\mathcal{G})$, is the maximum graph distance between any pair of vertices in $\mathcal{V}$, i.e.
\begin{equation}
diam(\mathcal{G}) = \max\{dist(v_{i}, v_{j}), v_{i}, v_{j}\in \mathcal{V}\}.
\end{equation}

This concept is useful to explain why spectral filters of ChebNet are exactly $K$-localized.

\section{Laplacian matrix}\label{sec:laplacian}
\subsection{Properties}
Weight on the graph is an associated numerical value assigned to each edge of a graph. A weighted graph is a graph associated with a weight to each of its edges while an unweighted graph is one without weights on its edges. The Laplacian matrix (unnormalized Laplacian or combinatorial Laplacian) for an unweighted graph is 
\begin{equation}
L = D - A \in \mathbb{R}^{n \times n}.
\end{equation}

Analogously, weighted graph is 
\begin{equation}
L = D - W \in \mathbb{R}^{n \times n},
\end{equation}
where $W$ is weighted adjacent matrix.

In graph theory, a regular graph is a graph in which each vertex has the same number of neighbors, i.e. each node has the same degree. $k$‑regular graph is a regular graph with vertices of degree $k$.

When $\mathcal{G}$ is $k$-regular, it is easy to see that
\begin{equation}
\widehat{L} = I - \frac{1}{k}A = \frac{1}{k}L,
\end{equation}
or 
\begin{equation}
\widehat{L} = I - \frac{1}{k}W = \frac{1}{k}L.
\end{equation}

In addition, the other Laplacian matrix, namely signless Laplacian,  denoted by $L_{s}$, is defined as $L_{s} = D + A$.

Eigen decomposition, also known as spectral decomposition, is a method to decompose a matrix into a product of matrices involving its eigenvalues and eigenvectors. Assuming basis of $\mathcal{L}$ is $U = (u_{1}, u_{2},\ldots,u_{n})$, $u_{i} \in \mathbb{R}, i = 1,2,\ldots,n$. Considering the Laplacian matrix is real symmetric matrix, the spectral decomposition of Laplacian matrix is 

\begin{equation}
\mathcal{L}  = U\Lambda U^{-1} = U\Lambda U^{T},
\label{sepctraldocom}
\end{equation}	
where
\begin{equation*}
\Lambda = \left[
\begin{matrix}
\lambda_{1}   &    &    &       \\
& \lambda_{2}   &    &       \\
&  & \ddots &  \\
&       &  & \lambda_{n}      \\
\end{matrix}
\right]=diag([\lambda_{1},\ldots,\lambda_{n}]) \in \mathbb{R}^{n \times n}.
\end{equation*}

Usually, the Laplacian matrix we referred is normalized Laplacian \cite[Section 1.3]{Butler2008}. It is easy to see that the Laplacian matrix $L$ associated with an undirected graph is positive semi-definite: Let $f = \{f_{1}, f_{2}, \ldots , f_{n}\}$ be an arbitrary vector, then
\begin{equation*}
\begin{aligned}
	f^{T}Lf & =f^TDf-f^TWf=\sum_{i=1}^{n}D_{i,i}f_{i}^2-\sum_{i,j}^{n}f_{i}f_{j}W_{i,j} \\ & =\frac{1}{2}(\sum_{i=1}^{n}D_{i,i}f_{i}^2 -  2\sum_{i=1}^{n}\sum_{j=1}^{n}f_{i}f_{j}W_{i,j}+\sum_{j=1}^{n}D_{j,j}f_{j}^2)
	\\ & =\frac{1}{2}(\sum_{i=1}^{n}\sum_{j=1}^{n}W_{i,j}(f_{i} - f_{j})^2) \ge 0
\end{aligned}
\end{equation*}

These are basic facts that simply follow from $L$'s symmetric and positive semi-definite properties:
\begin{itemize}
\item $L$ of order $n$ have $n$ linearly independent eigenvectors.
\item The eigenvectors corresponding to different eigenvalues of $L$ are orthogonal to each other, and the matrix formed by these orthogonal eigenvectors normalized to the unit norm is an orthogonal matrix.
\item The eigenvectors of $L$ can be taken as real vectors.
\item The eigenvalues of $L$ are nonnegative.
\end{itemize}

\textbf{Laplacian operator:}
The Laplacian matrix essentially is a Laplacian operator on a graph. To illustrate this concept, we introduce the incidence matrix. The incidence matrix is a matrix that reflect the relationship between vertices and edges. Suppose \textbf{the direction of each edge in the graph is fixed (but the direction can be set arbitrarily)}, let $f = (f_{1}, f_{2}, f_{3},\ldots, f_{n})^{T}$ denote signal vector associated with the vertices $(v_{1}, v_{2}, v_{3},\ldots, v_{n})$, the incidence matrix of a graph, denoted by $\nabla$, is a $|\mathcal{E}| \times |\mathcal{V}|$ matrix, the incidence matrix is defined as follows:
\begin{equation}
\nabla_{i,j} = \left\{
\begin{array}{rcl}
\nabla_{i,j} = -1& & \textrm{if $v_{j}$ is the initial vertex of edge $e_{i}$}\\
\nabla_{i,j} = 1 & & \textrm{if $v_{j}$ is the terminal vertex of edge of $e_{i}$}\\
\nabla_{i,j} = 0 & & \textrm{if $v_{j}$ is not in $e_{i}$}\\
\end{array} \right.
\end{equation}

\begin{figure}[ht]
	\centering
	\subfigure{\includegraphics[width=0.8\textwidth]{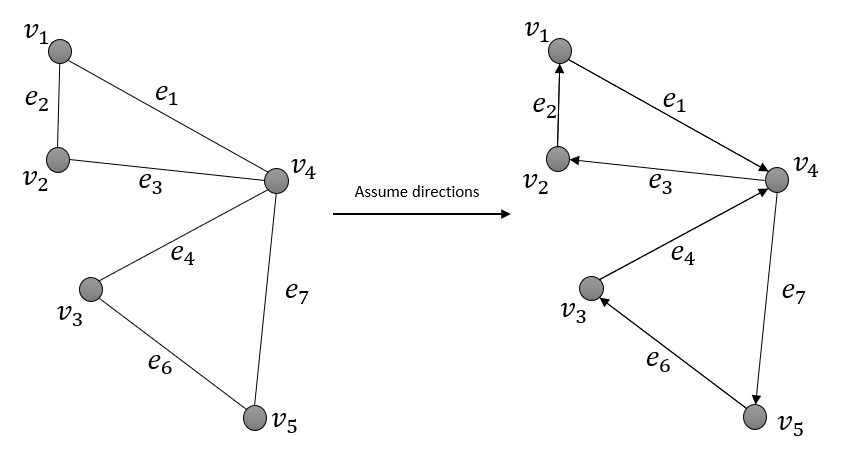}} 
	\caption{Graph}
	\label{gsp1}
\end{figure}

The mapping $f \longrightarrow \nabla f$ is known as the co-boundary mapping of the graph,  we take an example from the graph as shown in \figurename~\ref{gsp1}, we arrange arbitrary directions to the edges as the figure in right shows. We have
\begin{equation*}
\nabla = \left[
\begin{matrix}
-1 & 0 & 0 & 1 & 0 \\
1 & -1 & 0 & 0 & 0 \\
0 & 1 & 0 & -1 & 0 \\
0 & 0 & -1 & 1 & 0 \\
0 & 0 & 0 & -1 & 1 \\
0 & 0 & 1 & 0 & -1 \\
\end{matrix}\right].
\end{equation*}

Accordingly,
$\nabla\left[\begin{matrix}
f_1\\
f_2\\
f_3\\
f_4\\
f_5\\
\end{matrix}\right] 
=\left[
\begin{matrix}
-1 & 0 & 0 & 1 & 0 \\
1 & -1 & 0 & 0 & 0 \\
0 & 1 & 0 & -1 & 0 \\
0 & 0 & -1 & 1 & 0 \\
0 & 0 & 0 & -1 & 1 \\
0 & 0 & 1 & 0 & -1 \\
\end{matrix}\right]
\left[\begin{matrix}
f_1\\
f_2\\
f_3\\
f_4\\
f_5\\
\end{matrix}\right]
=\left[\begin{matrix}
f_4 - f_1\\
f_1 - f_2\\
f_2 - f_4\\
f_4 - f_3\\
f_5 - f_4\\
f_3 - f_5\\
\end{matrix}\right]
$ 

Therefore, $(\nabla f)(e_{i,j}) $ is given by
\begin{equation}
(\nabla f)(e_{i,j}) = f_{j} - f_{i}
\end{equation}
where $e_{i,j}$ denote the edge connecting node $i$ and node $j$.

Furthermore, 
\begin{equation*}
\nabla^{T} (\nabla f) = \left[
\begin{matrix}
-1 & 1 & 0 & 0 & 0 & 0\\
0 & -1 & 1 & 0 & 0 & 0\\
0 & 0 & 0 & -1 & 0 & 1\\
1 & 0 & -1 & 1 & -1 & 0\\
0 & 0 & 0 & 0 & 1 & -1\\
\end{matrix}\right] \left[\begin{matrix}
f_4 - f_1\\
f_1 - f_2\\
f_2 - f_4\\
f_4 - f_3\\
f_5 - f_4\\
f_3 - f_5\\
\end{matrix}\right] = \left[\begin{matrix}
2f_1 - f_2 -f_4\\
2f_2 - f_1 - f_4\\
2f_3 - f_4 - f_5\\
4f_4 - f_1 - f_2 - f_3 - f_5\\
2f_5 - f_3 - f_4\\
\end{matrix}\right]
\end{equation*}

Thus, the Laplacian matrix $L$ operating on $g$ would become
\begin{equation*}
\begin{aligned}
Lf = (D - A)f &= \left[
\begin{matrix}
2 & -1 & 0 & -1 & 0 \\
-1 & 2 & 0 & -1 & 0 \\
0 & 0 & 2 & -1 & -1 \\
-1 & -1 & -1 & 4 & -1 \\
0 & 0 & -1 & -1 & 2 \\
\end{matrix}\right]f = \left[\begin{matrix}
2f_1 - f_2 -f_4\\
2f_2 - f_1 - f_4\\
2f_3 - f_4 - f_5\\
4f_4 - f_1 - f_2 - f_3 - f_5\\
2f_5 - f_3 - f_4\\
\end{matrix}\right] \\&= \nabla^{T} (\nabla f)
\end{aligned}
\end{equation*}

Therefore, for any undirected graph, 
\begin{equation}
L f = \nabla^{T} (\nabla f) .
\end{equation} 

Particularly, for an $n$-dimensional Euclidean space, the Laplacian operator can be considered as a second-order differential operator

Analogously, consider undirected weighted graphs $\mathcal{G}$, each edge $e_{i,j}$ is weighted by $w_{i,j} > 0$, the Laplace operator on the graph can be defined as
\begin{equation}
(L f)_{i} = \sum_{j=1}^{n}W_{i,j}(f_{i} - f_{j}),
\end{equation} 
where $W_{i,j} = 0$ if $e_{i,j} \in \mathcal{E}$.

Also, 
\begin{equation*}
(L f)_{i} = \sum_{j}^{n}W_{i,j}(f_{i} - f_{j}) = D_{ii}f_{i} - \sum_{j}^{n} W_{i,j}f_{j} = (Df - Wf)_{i} = (L f)_{i}.
\end{equation*} 

For any $i$ holds, then it can be general form:
\begin{equation}
(D - W)f = Lf.
\end{equation}

As a quadratic form,
\begin{equation}
f^{T}Lf = \frac{1}{2}\sum_{e_{i,j}}W_{i,j}(f_{i} - f_{j})^{2}.
\end{equation}

Therefore, graph Laplacian matrix $L$, intrinsically as a Laplacian operator, make the centre node subtracts the surrounding nodes in turn, multiplying the corresponding link weights at the same time, and then sums them.

\subsection{Normalization}
In some cases of practical application, Laplacian matrix requires some kinds of normalization to ensure the algorithm convergence. Normalization of Laplacian matrix $\mathcal{L}$ include:
\begin{itemize}
\item Random-walk normalization:
\begin{equation}
\mathcal{L} = I - D^{-1}A
\end{equation}
\item  Symmetric normalization:
\begin{equation}
\mathcal{L} = D^{-1/2}LD^{-1/2} = I - D^{-1/2}WD^{-1/2},
\end{equation}
with the convention $D^{-1}_{i,i}=0$ for $d_{i}=0$, particularly, node $i$ is an isolated vertex if $d_{i}=0$. A graph is said to be non-trivial if it contains at least one edge \cite[Section 1.2]{Chung:1997}.  The normalized Laplacian has eigenvalues always lying in the range between 0 and 2 inclusive as demonstrated by Chung \cite[Section 1.3]{Chung:1997}.
\end{itemize}
In the following, we will use symmetric normalized Laplacian matrix as default, unless unless otherwise stated. 

\subsection{Eigenvalues discussion}
The spectral GCN analysis relies on spectral graph theory, which studies the properties of graphs via the eigenvalues and its corresponding eigenvectors associated with the graph adjacency matrix and the graph Laplacian matrix and its variants. Since graph convolutional operator is defined on eigenvalues of the Laplacian matrix, being familiar with the properties of graph-related properties (as the facts listed below) is greatly helpful for the research on spectral GCN. 

The eigenvalues of the Laplacian matrix can be inferred from the properties of the graph, for example:
\begin{lemma}[\protect{\cite[Section 1.3]{Chung:1997}}]
	The number of zero eigenvalues of the Laplacian (i.e. the multiplicity of 0 as an eigenvalue) is the number of connected components \footnote{A connected component or simply component of an undirected graph is a subgraph in which each pair of nodes is connected with each other via a path.} of $\mathcal{G}$. In particular $\mathcal{G}$ is connected if and only if $\lambda_{2} > 0$ ($\lambda_{2}$ is the second smallest eigenvalue, some references use $\lambda_{1}$). The multiplicity of 2 as an eigenvalue is the number of bipartite connected components of $\mathcal{G}$ with at least two vertices. 
\end{lemma}

\begin{theorem}\label{thm:eigenrange}
	The matrix $\mathcal{L}$ is positive semi-definite and satisfies: All eigenvalues lie in the interval $[0,2]$.
\end{theorem}

Theorem~\ref{thm:eigenrange} gives us the range of the eigenvalues of the normalized Laplacian matrix $\mathcal{L}$. For a stable graph filter of GCNs, it requires the absolute eigenvalues of $\mathcal{L}$ to be bounded by 1, thus various scaling methods are introduced as follows.

\section{Discrete Signal Processing on Graphs}\label{sec:discretespg}

\textbf{Graph filters}, generally in discrete signal processing(DSP), is a system $H(\cdot)$ that takes a graph signal $f$ as an input, processes it, and produces another graph signal $\tilde{f} = H(f)$ as an output \cite{10.1109/TSP.2014.2321121}, \cite{6638850}. In discrete signal processing on graphs (DSPG) \cite{6638850}, an equivalent concept of filters for the processing of graph signals. Given graph signals $f$ indexed by a graph , the fundamental building block for graph filters on $G$ is a graph shift that replaces each signal coefficient $f_{i}$ indexed by node $n$ with a linear combination of coefficients at other nodes weighted proportionally to the degree of their relation \cite{6638850}:
\begin{equation}
\tilde{f}_{i} = \sum_{m=1}^{n} W_{i,m}f_{m} \Leftrightarrow \tilde{f} = Wf,
\end{equation}
where W is the graph shift or a weighted adjacency matrix.

According to \cite[Theorem~1]{10.1109/TSP.2013.2238935}, any linear, shift-invariant graph filter is necessarily a matrix polynomial in the (weighted) adjacency matrix $W$ of the form 
\begin{equation}
h(W) = h_{0}I + h_{1}W + \ldots + h_{L}W^{L}.
\label{polyfilter}
\end{equation}
The output of the filter \eqref{polyfilter} is the signal
\begin{equation}
\tilde{f} = H(f) = h(W)f,
\end{equation}
where $h_{l} \in \mathbb{C}$ are possible coefficients. In addition, $L \le n$, which means any graph filter  \refeq{polyfilter} can be represented by at most $n$ coefficients. Also, if graph filter \refeq{polyfilter} is invertible, matrix $h(A)$ is non-singular, its inverse also is a matrix polynomial in $W$ of the form \refeq{polyfilter}, namely $g(W) = h(W)^{-1}$. 

Generally, a Fourier transform is a uniform to the expansion of a signal using basis elements that are invariant to filtering. And the basis can be the eigenbasis of the $W$ as \refeq{sepctraldocom} or the Jordan eigenbasis of $W$ if the complete eigenbasis does not exist.


\textbf{Graph Fourier transform} is analogous to classical Fourier transform, similarly, the eigenvalues could represent graph frequencies and form the spectrum of the graph,  eigenvectors denote frequency components which serve the work as the graph Fourier basis \cite{6638850}, \cite{7208894}. 

Let graph Fourier basis $U = (u_{1}, u_{2},\ldots,u_{n})$, $ u_{i} \in \mathbb{R}, i = 1,2,\ldots,n$ from Laplacian matrix $\mathcal{L}$ as \refeq{sepctraldocom}.  Nodes' signal $f = (f_{1}, f_{2}, f_{3},\ldots, f_{n})^{T}$, after graph Fourier Transform, signal become $\hat{f} = (\hat{f}(\lambda_{1}), \hat{f}(\lambda_{2}), \hat{f}(\lambda_{3}),\ldots, \hat{f}(\lambda_{n}))^{T}$, the graph Fourier transform is
\begin{equation}
\hat{f} = U^{T}f.
\end{equation}

Correspondingly, inverse Graph Fourier transform is
\begin{equation}
f = U\hat{f}.
\end{equation}

Therefore, taking Laplace's eigenvector as the basis function, any signal on the graph can be
\begin{equation}
f = \hat{f}(\lambda_{1}) u_{1} + \hat{f}(\lambda_{2}) u_{2} + \ldots + \hat{f}(\lambda_{n}) u_{n} = \sum_{i=1}^{n}\hat{f}(\lambda_{i})u_{i},
\end{equation}
$u_{i}$ is the column vector of orthogonal matrix from spectral decomposition from $\mathcal{L} = U\Lambda U^{T}$.

In fact, that is analogous to the principle of Discrete Fourier Transform(DFT)
\begin{equation}
X_{2\pi}(k)=\sum_{n=-\infty}^{\infty}x_{n}e^{-ikn}.
\end{equation}

\section{Spectral graph convolution}\label{sec:sgc}
\subsection{Overview}

The principal of convolutional neural network is beyond the discussion of this paper, we refer the readers to \cite{DLYANN} for a fundamental understanding. In the following, we will define graph filter, which is an convolution operator on graph in the fourier domain, as well as other associated concepts, which leads to various classic GCN models. 

In the Fourier domain, the convolution operator on graph $\cdot_{G}$ is defined as
\begin{equation}
g(\cdot_{G})f = \mathcal{F}^{-1}(\mathcal{F}(g) \odot \mathcal{F}(f)) = U(U^{T}g \odot U^{T}f) = U g_{\theta}(\Lambda)  U^{T}f = g_{\theta}(\mathcal{L})f.
\end{equation}
where $(\cdot_{G})$ is convolution operator defined on graph, $\odot$ is Hadamard product.

It follows that a signal $f$ is filtered by $g \in \mathbb{R}^n$, and denotes $g_{\theta}(\Lambda) = diag(U^{T}g)$ which the diagonal corresponds to spectral filter coefficients. 

For details, 
\begin{equation*}
\begin{aligned}
	g_{\theta}(\cdot_{G})f & = g_{\theta}(\mathcal{L})f = g_{\theta}(U\Lambda U^{T})f = Ug_{\theta}(\Lambda) U^{T}f \\& =U\left[
	\begin{matrix}
	\hat{g}(\lambda_{1})   &    &    &       \\
	& \hat{g}(\lambda_{2})   &    &       \\
	&  & \ddots &  \\
	&       &  & \hat{g}(\lambda_{n})      \\
	\end{matrix}
	\right]U^{T}f \\& = U\left[
	\begin{matrix}
	\hat{g}(\lambda_{1})   &    &    &       \\
	& \hat{g}(\lambda_{2})   &    &       \\
	&  & \ddots &  \\
	&       &  & \hat{g}(\lambda_{n})      \\
	\end{matrix}
	\right]\hat{f} \\& = U\left[
	\begin{matrix}
	\hat{g}(\lambda_{1})   &    &    &       \\
	& \hat{g}(\lambda_{2})   &    &       \\
	&  & \ddots &  \\
	&       &  & \hat{g}(\lambda_{n})      \\
	\end{matrix}
	\right]\left[
	\begin{matrix}
	\hat{f}(\lambda_{1})\\
	\hat{f}(\lambda_{2})\\
	\ldots\\
	\hat{f}(\lambda_{n})\\
	\end{matrix}\right] \\& = U\left[
	\begin{matrix}
	\hat{g}(\lambda_{1})\\
	\hat{g}(\lambda_{2})\\
	\ldots\\
	\hat{g}(\lambda_{n})\\
	\end{matrix}\right] \odot \left[
	\begin{matrix}
	\hat{f}(\lambda_{1})\\
	\hat{f}(\lambda_{2})\\
	\ldots\\
	\hat{f}(\lambda_{n})\\
	\end{matrix}\right]. 
\end{aligned}
\end{equation*}

Spectral-based GCN all follow this definition of $Ug_{\theta}(\Lambda) U^{T}f$, the main difference between different version of Spectral-based GCN lies in the choice of the filter $g_{\theta}(\Lambda)$ \cite{Wu2021ACS}.

\subsection{Spectral CNN}

Bruna et al. propose the first spectral convolutional neural network \cite{ae482107de73461787258f805cf8f4ed}.
A graph can be associated with node signal $f \in \mathbb{R}^{n \times C_{k}}$ is a feature matrix with $f_i \in \mathbb{R}^{C_{k}}$ representing the feature vector of node $i$. A construction where each layer $k=1,\ldots,K$ transforms an input vector $f^{(k)}$ of size $n \times C_{k}$ into an output $f^{(k + 1)}$ of size $n \times C_{k+1}$.

\begin{equation}
f_{j}^{(k + 1)}=\sigma(U\sum_{i=1}^{C_{k}}{g_{\theta}}_{i,j}^{(k)}U^{T}f_{i}^{(k)}) = \sigma(U\sum_{i=1}^{C_{k}}{g_{\theta}}_{i,j}^{(k)}\hat{f}_{i}^{(k)}),
\end{equation}
where ${g_{\theta}}_{i,j}^{(k)}, i=1,\ldots,n; j=1,\ldots,C_{k}$ is a diagonal matrix with trainable parameters $\theta_{m}^{(k)}, m \in (1,n)$, $\sigma$ is activation function. ${g_{\theta}}_{i,j}^{(k)}$ is given by
\begin{equation*}
{g_{\theta}}_{i,j}^{(k)} = \left[
\begin{matrix}
\theta_{1}^{(k)}   &    &    &       \\
& \theta_{2}^{(k)}   &    &       \\
&  & \ddots &  \\
&       &  & \theta_{n}^{(k)}    \\
\end{matrix}
\right].
\end{equation*}

\subsection{ChebNet}

ChebNet \cite{10.5555/3157382.3157527} uses Chebyshev polynomials instead of convolutions in spectral domain. Furthermore, it was demonstrated that that $g_{\theta}(\Lambda)$ can be approximated by a truncated expansion in terms of Chebyshev polynomials \cite{Hammond:131283}.
\begin{equation}
T_{n+1}(x)=2xT_{n}(x) - T_{n-1}(x), n \in \mathbb{N}^{+},
\label{chebnetpoly}
\end{equation}
where $T_{0}(x) = 1, T_{1}=x$. Here, we make $\widetilde{\Lambda} = \frac{2\Lambda}{\lambda_{max}} - I_{n} \in [-1,1]$, $\lambda_{max}$ is the biggest eigenvalue from $\mathcal{L}$
\begin{equation}
g_{\theta}(\Lambda) = \sum_{k=0}^{K-1}\theta_{k}T_{k}(\widetilde{\Lambda}),
\end{equation}
where the parameter $\theta \in \mathbb{R}^{K}$.

The filtering operator can also be written as 
\begin{equation}
g_{\theta}(\mathcal{L})f =  \sum_{k=0}^{K-1} \theta_{k}T_{k}(\tilde{\mathcal{L}})f,
\label{eqchebnet}
\end{equation}
where $T_{k}(\tilde{\mathcal{L}}) \in \mathbb{R}^{n \times n}$ is the Chebyshev polinomial of order $k$ evaluated at the scaled Laplacian $\tilde{\mathcal{L}} = 2 \mathcal{L}/\lambda_{max} - I_{n}$. Accordingly, spectral filters represented by $K^{th}$-order polynomials of the Laplacian are exactly $K$-localized, i.e. it depends only on nodes that are at maximum $K$ steps away from the central node \cite{10.5555/3157382.3157527}, \cite[Lemma 5.2]{Hammond:131283}. 

\begin{lemma}[\protect{\cite[Lemma 5.2]{Hammond:131283}}]
	Let $\mathcal{G}$ be a weighted graph, with adjacency matrix $A$. Let $B$ equal the adjacency matrix of the binarized graph, i.e. $B_{m,n} = 0$ if $A_{m,n} = 0$, and $B_{m,n} = 1$ if $A_{m,n} > 0$. Let $\tilde{B}$ be the adjacency matrix with unit loops added on every vertex, e.g. $\tilde{B}_{m,n} = B_{m,n}$ for $m \neq n$ and $\tilde{B}_{m,n}=1$ for $m=n$.
	
	Then for each $s > 0$, $(B^{s})_{m,n}$ equals the number of paths of length $s$ connecting $m$ and $n$, and $(\tilde{B}^{s})_{m,n}$ equals the numebr of all paths of length $r \le s$ connecting $m$ and $n$.
\end{lemma}
The Lemma can be used to demonstrate that matrix elements of low powers of the graph Laplacian corresponding to sufficiently separated vertices must be zero. Therefore, $dist(v_{i}, v_{j}) > K$ implies $(\mathcal{L}^K)_{i,j} = 0$, and the spectral filters of ChebNet are exactly $K$-localized.

Accordingly,
\begin{equation*}
\begin{aligned}
g_{\theta}(\Lambda) &= \left[
	\begin{matrix}
	\hat{g}(\lambda_{1})   &    &    &       \\
	& \hat{g}(\lambda_{2})   &    &       \\
	&  & \ddots &  \\
	&       &  & \hat{g}(\lambda_{n})      \\
	\end{matrix}
	\right] \\& = \left[
	\begin{matrix}
	\sum_{k=0}^{K-1}\theta_{k}T_{k}(\hat{\lambda_{1}})   &    &    &       \\
	& \sum_{k=0}^{K-1}\theta_{k}T_{k}(\hat{\lambda_{2}})   &    &       \\
	&  & \ddots &  \\
	&       &  & \sum_{k=0}^{K-1}\theta_{k}T_{k}(\hat{\lambda_{n}})      \\
	\end{matrix}
	\right],
\end{aligned}
\end{equation*}
where $\theta_{k}$ is a vector of Chebyshev coefficients, which is trainable parameter.

Furthermore, Equation \refeq{eqchebnet} can be deduced as following
\begin{equation}
\begin{aligned}
f(\cdot_{G}) g_{\theta} = & g_{\theta}(U\Lambda U^{T})f = U\sum_{k=0}^{K-1}\theta_{k}T_{k}(\widetilde{\Lambda})U^{T}f =\sum_{k=0}^{K-1}U\theta_{k}T_{k}(\widetilde{\Lambda})U^{T}f\\ &
=\sum_{k=0}^{K-1}U\theta_{k}(\sum_{c=0}^{k}\alpha_{kc}\widetilde{\Lambda}^{k})U^{T}f
=\sum_{k=0}^{K-1}\theta_{k}(\sum_{c=0}^{k}\alpha_{kc}U\widetilde{\Lambda}^{k}U^{T})f\\&
=\sum_{k=0}^{K}\theta_{k}(\sum_{c=0}^{k}\alpha_{kc}(U\widetilde{\Lambda}U^{T})^{k})f
=\sum_{k=0}^{K-1}\theta_{k}T_{k}(U\widetilde{\Lambda}U^{T})f \\ & =\sum_{k=0}^{K-1}\theta_{k}T_{k}(\widetilde{\mathcal{L}})f.
\end{aligned}
\end{equation}

After using Chebyshev polynomial instead of the convolution kernel of the spectral domain, ChebNet does not need the Laplace matrix is to be eigen-decomposed. The most time-consuming steps are omitted \cite{10.5555/3157382.3157527}.

\textbf{Comparison between Spectral CNN and ChebNet}

Assuming that $n$ is the number of nodes.
\begin{itemize}
	\item The parameter complexity of the SCNN model is very large, and the learning complexity is $O(n)$ \cite{ae482107de73461787258f805cf8f4ed}, \cite[Section 2.1]{10.5555/3157382.3157527}, which is easy to overfit when there are many nodes. When dealing with large-scale graph data which usually has more than millions of nodes, it will face great challenges.
	\item Computing the eigenvalue decomposition of the Laplace matrix is very time-consuming.
	\item The convolution kernel of ChebNet has only K learnable parameters($\theta_{k}$), and $K \ll n$, hence their learning complexity is $O(K)$, the complexity of learnable parameters is greatly reduced \cite[Section 2.1]{10.5555/3157382.3157527}.
	\item ChebNet does not need the Laplace matrix to be eigen-decomposed, instead it approximate $g_{\theta}(\mathcal{L})$ with a truncated expansion in term of Chebyshev polynomials $T_{k}(x)$ of $K^{th}$ order \cite[Section 2.1]{10.5555/3157382.3157527}.
\end{itemize}

\subsection{CayleyNets}
The paper \cite{CayleyNets2017} construct a family of complex filters that enjoy the advantages of Chebyshev filters while avoiding some of their drawbacks. A Cayley polynomial of order $r$ to be a real-valued function with complex coefficients.
\begin{equation}
g_{c,h}(\lambda) = c_{0} + 2Re\{\sum_{j=1}^{r}c_{j}(h\lambda - i)^{j}(h\lambda + i)^{-j}\},
\end{equation}
where $c=(c_{0,\ldots,c_{r}}) $ is a vector of one real coefficient and $r$ complex coefficients and $h > 0$ is the spectral zoom parameter. 

A Cayley filter $G$ is a spectral filter defined on real signals $f$ by 
\begin{equation}
g_{\theta}(\mathcal{L})f = g_{c,h}(\Lambda)f =  c_{0}f + 2Re\{\sum_{j=1}^{r}c_{j}(h\mathcal{L} - iI)^{j}(h\mathcal{L} + iI)^{-j}f\}.
\end{equation}
the parameters c and h is learnable, which are optimized during training.

\vspace{10pt}

The application of the filter $g_{\theta}(\mathcal{L})f$ can be performed without explicit expensive eigendecomposition of the Laplacian operator. The unit complex circle is denoted by $e^{i\mathbb{R}} = \{e^{i\theta}, \theta \in \mathbb{R}\}$.

\vspace{10pt}

The Cayley transform $C(x)=\frac{x-i}{x+i}$ is a smooth bijection between $\mathbb{R}$ and $e^{i\mathbb{R}} \setminus \{1\}$.

\vspace{10pt}

Correspondingly, by applying the Cayley transform to the scaled Laplacian $h\mathcal{L}$, we get the complex matrix 
\begin{equation}
C(h\mathcal{L})=(h\mathcal{L} - iI)(h\mathcal{L} - iI)^{-1}.
\end{equation}
which has its spectrum in $e^{i\mathbb{R}}$ and is thus unitary.

Since $z^{-1}=\overline{z}$ for $z \in e^{i\mathbb{R}}$, we have  $\overline{c_{j}C^{j}(h\mathcal{L})}=\overline{c_{j}}C^{-j}(h\mathcal{L})$ and given $2Re{z} = z + \overline{z}$, any Cayley filter can be written as a conjugate-even Laurent polynomial. 
\begin{equation}
g_{\theta} = c_{0}I + 2Re\{\sum_{j=1}^{r}c_{j}(h\mathcal{L} - i)^{j}(h\mathcal{L} + i)^{-j}f\}.
\end{equation}

\textit{proof}:
\begin{equation*}
\begin{aligned}
g_{\theta}(\mathcal{L}) &= c_{0}I + \sum_{j=1}^{r}[c_{j}C^{j}(h\Delta) + \overline{c_{j}}C^{-j}(h\Delta)] \\&=c_{0}I + \sum_{j=1}^{r}[c_{j}C^{j}(h\Delta) + \overline{c_{j}C^{j}(h\Delta)}] \\&= c_{0}I + 2Re\{\sum_{j=1}^{r}c_{j}(h\mathcal{L} - iI)^{j}(h\mathcal{L} + iI)^{-j}f\}.
\end{aligned}
\end{equation*}

Since the spectrum of $C(h\Delta)$ is in $e^{i\mathbb{R}}$, the operator $C^{j}(h\Delta)$ can be thought of as a multiplication by a pure harmonic in the frequency domain $e^{i\mathbb{R}}$ for any integer power $j$, 
\begin{equation}
C^{j}(h \mathcal{L}) = U diag([C(h\lambda_{1})]^{j},\ldots,[C(h\lambda_{n})]^{j})U^{T},
\end{equation}
where $C^{j}(h \mathcal{L})$ it is a (real-valued) trigonometric polynomial, and $g_{\theta}(\Lambda)$ is conjugate-even. 

Hence, a Cayley filter $g_{\theta}$ can be seen as a multiplication by a finite Fourier expansion in the frequency domain $e^{i\mathbb{R}}$. While preserving spatial locality, ChebNet can be considered as a special case of CayleyNet \cite{CayleyNets2017}, \cite{Wu2021ACS}.

\subsection{GCN}
GCN \cite{Kipf:2016tc} can be regarded as a further simplification of ChebNet. To reduce the computational complexity, only the first order Chebyshev polynomials are considered, consequently each convolution kernel has only one trainable parameter \cite{Kipf:2016tc}. Combining with \eqref{chebnetpoly}, we have
\begin{equation}
g_{\theta}(\Lambda)=\sum_{k=0}^{1}\theta_{k}T_{k}(\widetilde{\Lambda}).
\end{equation}

Hence,
\begin{equation}
g_{\theta}(\Lambda) = \left[
\begin{matrix}
\sum_{k=0}^{1}\theta_{k}T_{k}(\hat{\lambda_{1}})   &    &    &       \\
& \sum_{k=0}^{1}\theta_{k}T_{k}(\hat{\lambda_{2}})   &    &       \\
&  & \ddots &  \\
&       &  & \sum_{k=0}^{1}\theta_{k}T_{k}(\hat{\lambda_{n}})      \\
\end{matrix}
\right].
\end{equation}

\vspace{10pt}

In this linear formulation of a GCN we further approximate $\lambda_{max} \approx 2$. Under such approximations, this can simplifies to:
\begin{equation}
\widetilde{\mathcal{L}}=\frac{2}{\lambda_{max}}\mathcal{L} - I_{n} = \mathcal{L} - I_{n},
\end{equation}
where $\mathcal{L}$ is normalized graph Laplacian $\mathcal{L} = I - D^{-\frac{1}{2}}AD^{-\frac{1}{2}}$.

Then,
\begin{equation}
f(\cdot_{G}) g = \sum_{k=0}^{1}\theta_{k} T_{k}(\widetilde{\mathcal{L}})f =\theta_{0}T_{0}(\widetilde{\mathcal{L}})f+\theta_{1}T_{1}(\widetilde{\mathcal{L}})f,
\end{equation}
where $A$ is an adjacency matrix of the graph.
\vspace{10pt}

Accordingly, 
\begin{equation}
f(\cdot_{G})g = (\theta_{0}+\theta_{1}(\mathcal{L}-I_{n}))f=(\theta_{0}-\theta_{1}(D^{-\frac{1}{2}}AD^{\frac{1}{2}}))f.
\end{equation}

Furthermore, to reduce the number of trainable parameters——each kernel has only one trainable parameter, we set $\theta_{0}=-\theta_{1}=\theta$, then we have
\begin{equation*}
f(\cdot_{G}) g \approx (\theta_{0}+\theta_{1}(\mathcal{L} -I_{n}))f = \theta_{0} - \theta_{1}D^{-\frac{1}{2}}AD^{-\frac{1}{2}}=(\theta(D^{-\frac{1}{2}}AD^{-\frac{1}{2}} + I_{n}))f,
\end{equation*}
where $D^{-\frac{1}{2}}AD^{-\frac{1}{2}} + I_{n}$ now has eigenvalues in the range [0, 2]. Then, only one parameter in convolution kernel can be learned. The number of parameters is greatly reduced, which can reduce the number of parameters to prevent overfitting.

However, repeated application of this operator can therefore lead to numerical instabilities and exploding or vanishing gradients. To alleviate this problem, the following re-normalization trick is introduced.

We add self-loop to $A$, 
\begin{equation}
\widetilde{A} = A + I_{n}.
\end{equation}

Correspondingly,
\begin{equation}
\widetilde{D}_{i,i} = \sum_{j=1}^{n}\widetilde{A}_{i,j}.
\end{equation}

Finally, 
\begin{equation}
f(\cdot_{G}) g=\theta \widetilde{D}^{-\frac{1}{2}}\widetilde{A}\widetilde{D}^{-\frac{1}{2}}f.
\end{equation}

Usually, we write $\theta$ as $W, \hat{A} = \widetilde{D}^{-\frac{1}{2}}\widetilde{A}\widetilde{D}^{-\frac{1}{2}}$, then we have $f(\cdot_{G}) g = \hat{A} f W$.

Here make an illustration of example(applied on Cora dataset, a node level task), consider a two-layer GCN for semi-supervised node classification on a graph, $f \in \mathbb{R}^{n \times C}$ is $n$ nodes with $C$ input channels
\begin{equation}
Z = softmax(\hat{A}  ReLU(\hat{A}fW^{<0>}) W^{<1>}),
\end{equation}

$W^{<0>} \in \mathbb{R}^{C \times H}$ is an input-to-hidden weight matrix for a hidden layer with H feature maps. $W^{<1>} \in \mathbb{R}^{H \times F}$ is a hidden-to-output weight matrix, $F$ is the dimension of feature maps in the output layer.

For calculation of loss function, need to evaluate the cross-entropy error over all labeled examples
\begin{equation}
Loss=-\sum_{l \in \gamma_{L}} \sum_{f=1}^{F}Y_{l,f}ln Z_{l,f},
\end{equation}
where $\gamma_{L}$ is the set of node indices that have labels, labels are denoted by $Y_{i}$. we then can use Stochastic Gradient descent as optimizer to finish the process of training.

\section{Accelerated filtering using Lanczos method}\label{sec:lanczos}
Given graph $\mathcal{G}$ and its corresponding Laplacian matrix $\mathcal{L}$, a non-zero vector $f \in \mathbb{R}^{n}$, we apply Lanczos algorithm \cite[Section 10.2]{GoluVanl96} as shown in Algorithm \ref{algorithmlanczos} to compute an orthonormal basis $V_{M}=[v_{1},\ldots, v_{M}]$ of the Krylov subspace $K_{M}(\mathcal{L},f)=span\{f, \mathcal{L}f,\ldots,\mathcal{L}^{M-1}f\}$.

Lanczos algorithm can form a symmetric tridiagonal matrix $H_{M} \in \mathbb{R}^{M \times M}$ 
\begin{equation}
V_{M}^{\ast}\mathcal{L}V_{M}=H_{M}= \left[
\begin{matrix}
\alpha_{1}  &  \beta_{2}   &    &   &   \\
\beta_{2} &  \alpha_{2}   &  \beta_{3}  &   & \\
& \beta_{3}  &  \alpha_{3}  &     \ddots  &\\
&   & \ddots & \ddots & \beta_{M} \\
&     &   & \beta_{M} & \alpha_{M}     \\
\end{matrix}
\right].
\end{equation}

\vspace{10pt}

\begin{algorithm}[H]
	\SetAlgoLined
	\large
	\KwIn{Symmetric matrix $\mathcal{L} \in \mathbb{R}^{n \times n}$, vector $f \ne 0$, $M \in \mathbb{N}$}
	\KwResult{$V_{M}=[v_{1}, \ldots, v_{M}]$ with orthonormal columns, scalars $\alpha_{1},\ldots,\alpha_{M},\beta_{2},\ldots,\beta_{M} \in \mathbb{R}$.}
	$v_{1}$ $\gets$ $f/\| f\|_{2} $\;
	\For {$j:=1$  \KwTo  $M$}{ 
		$w$ = $\mathcal{L}v_{j}$\;
		$\alpha_{j}$ = $v_{j}^{\ast}w$\;
		$\tilde{v}_{j+1}$ = $w$ - $v_{j}\alpha_{j}$\;
		\If{$j$ $>$ 1}{
			$\tilde{v}_{j+1}$ $\gets$ $\tilde{v}_{j+1}$ - $v_{j-1}\beta_{j-1}$\;
		}
		$\beta_{j}$ = $\| \tilde{v}_{j+1} \|_{2}$\;
		$\tilde{v}_{j+1}$ = $\tilde{v}_{j+1} / \beta_{j}$
	}
	\label{algorithmlanczos}
	\caption{\large Lanczos method}
\end{algorithm}

\vspace{20pt}

The approximation to $g_{\theta}(\mathcal{L})f$ is given by \cite{doi:10.1137/0913071}, \cite{Susnjara2015}
\begin{equation}
g_{\theta}(\mathcal{L})f \approx \| f \| _{2} V_{M}g_{\theta}(H_{M})e_{1}:=g_{M},
\end{equation}
where $e_{1} \in \mathbb{R}^{M}$ is the first unit vector. Because eigenvalue interlacing \footnote{Consider two sequences if real numbers: $\lambda_{1} \ge \lambda_{2} \ldots \ge \lambda_{n}$, and $\mu_{1} \ge \mu_{2} \ldots \ge \mu_{n}$ with $m < n$. The second sequence is said to interlace the first one whenever $\lambda_{i} \ge \mu_{i} \ge \lambda_{n-m+i}$ for $i = 1, \ldots, m$.} \cite{Haemers95interlacingeigenvalues}, the eigenvalues of $H_M$ are contained in the interval $[0, \lambda_{max}]$ and hence the expression $g_{\theta}(H_{M})$ is well-defined \cite{Susnjara2015}. Particularly, $M \ll n$, the computational cost by evaluating $g_{\theta}(H_{M})$ is inexpensive.

\begin{theorem}[\protect{\cite[Corollary 3.4]{NMOPS}}]
	Let $\mathcal{L} \in \mathbb{R}^{n \times n}$ be symmetric with eigenvalues contained in the interval [0, $\lambda_{max}$] and let $g_{\theta}:[0,\lambda_{max}] \rightarrow \mathcal{R}$ be continuous. Then
	\begin{equation}
	|| g_{\theta}(\mathcal{L}f - g_{M})\|_{2} \le 2 \| f\|_{2} \cdot \min_{p \in \mathcal{P}_{M-1}}\max_{z \in [0, \lambda_{max}]}|g_{\theta}(z) - p(z)|,
	\end{equation}
	where $\mathcal{P}_{M-1}$ denotes all polynomials of degree at most $M-1$.
	\label{theoremnmops}
\end{theorem}

According to Theorem \ref{theoremnmops}  the error is bounded by the best polynomial approximation \cite[Theorem 2.4.1]{Phillips2003InterpolationAA} of $g_{\theta}$ on $[0, \lambda_{max}]$. The paper \cite{Susnjara2015} demonstrates that{---}up to a multiple of two{---}the Lanczos-based approximation $g_{\theta}$ can be expected to provide at least the same accuracy. In addition, the Lanczos-based approximation can sometimes be expected to perform much better because of its ability to adapt to the eigenvalues of $\mathcal{L}$ \cite{Susnjara2015}, which can be well-understand for Krylov subspace approximations to solutions of linear systems \cite[Section 3.1]{doi:10.1137/1.9781611970937}.


\nocite{*}
\bibliographystyle{abbrv}
\bibliography{main}

\end{document}